\documentclass[11pt]{article}
\usepackage{amssymb,subfigure,amscd,amsthm,amsmath}

\usepackage{graphicx}

\usepackage[all,2cell]{xy}
\UseAllTwocells

\title{Investigating The Algebraic Structure of Dihomotopy Types}
\author{Philippe Gaucher\\
{\small Institut de Recherche Math\'ematique Avanc\'ee}\\
{\small ULP et  CNRS}\\
{\small 7 rue Ren\'e Descartes}\\
{\small 67084 Strasbourg}\\
{\small France}\\
{\small gaucher@math.u-strasbg.fr}}

\newcommand{\C}{\mathcal{C}}
\newcommand{\D}{\mathcal{D}}
\newcommand{\Z}{\mathbb{Z}}
\newcommand{\N}{\mathbb{N}}
\newcommand{\R}{\mathbb{R}}
\newcommand{\de}{\partial}
\newcommand{\p}\times
\renewcommand{\vec}{\overrightarrow}
\renewcommand{\P}{\mathbb{P}}

\newcommand{\Vt}{\mathcal{V}}

\newcommand{\be}{\begin{equation}}
\newcommand{\ee}{\end{equation}}
\newcommand{\bea}{\begin{eqnarray}}
\newcommand{\eea}{\end{eqnarray}}
\newcommand{\beas}{\begin{eqnarray*}}
\newcommand{\eeas}{\end{eqnarray*}}

\newtheorem{thm}{Theorem}[section]
\newtheorem{prop}[thm]{Proposition}

\newtheorem{question}[thm]{Question}

\newtheorem{rem}[thm]{Remarque}

\newtheorem{defn}{Definition}[section]
\newtheorem{propdef}[thm]{Proposition and Definition}
\newtheorem{philo}[thm]{Philosophy}

\newcommand{\bd}{\begin{defn}}
\newcommand{\ed}{\end{defn}}
\newcommand{\bcd}{\begin{defn}}
\newcommand{\ecd}{\end{defn}}
\newcommand{\bex}{\begin{exmp}}
\newcommand{\eex}{\end{exmp}}
\newcommand{\bp}{\begin{prop}}
\newcommand{\ep}{\end{prop}}
\newcommand{\bth}{\begin{thm}}
\renewcommand{\eth}{\end{thm}}
\newcommand{\br}{\begin{rem}}
\newcommand{\er}{\end{rem}}
\newcommand{\bpf}{\begin{proof}}
\newcommand{\epf}{\end{proof}}

\newcommand{\fl}[1]{\ar@{->}[l]_{#1}}
\newcommand{\fr}[1]{\ar@{->}[r]^{#1}}
\newcommand{\fd}[1]{\ar@{->}[d]_{#1}}
\newcommand{\fu}[1]{\ar@{->}[u]|{#1}}
\newcommand{\f}[2]{\ar@{->}[#1]|{#2}}
\newcommand{\ff}[2]{\ar@2{->}[#1]|{#2}}
\newcommand{\frr}[1]{\ar@{->}[rr]^{#1}}

\newcommand{\ei}{\underline{\iota}}
\newcommand{\ef}{\underline{\sigma}}

\newcommand{\CW}{{\mathbf{CW}}}
\newcommand{\diCW}{{\mathbf{glCW}}}
\newcommand{\lpohaus}{{\mathbf{LPoHaus}}}

\newcommand{\haus}{{\mathbf{Haus}}}
\renewcommand{\top}{{\mathbf{Top}}}
\newcommand{\potop}{{\mathbf{PoTop}}}

\newcommand{\Ho}{{\mathbf{Ho}}}
\newcommand{\iso}{\simeq}

\newcommand{\vI}{\vec{I}}

\date{July 2001}

\begin{document}

\maketitle

\begin{abstract}
  This presentation is the sequel of a paper published in GETCO'00
  proceedings where a research program to construct an appropriate
  algebraic setting for the study of deformations of higher
  dimensional automata was sketched. This paper focuses
  precisely on detailing some of its aspects. The main idea is that
  the category of homotopy types can be embedded in a new category of
  dihomotopy types, the embedding being realized by the Globe functor.
  In this latter category, isomorphism classes of objects are exactly
  higher dimensional automata up to deformations leaving invariant
  their computer scientific properties as presence or not of deadlocks
  (or everything similar or related). Some hints to study the
  algebraic structure of dihomotopy types are given, in particular a
  rule to decide whether a statement/notion concerning dihomotopy
  types is or not the lifting of another statement/notion concerning
  homotopy types. This rule does not enable to guess what is the
  lifting of a given notion/statement, it only enables to make the
  verification, once the lifting has been found.
\end{abstract}

\tableofcontents

\section{Introduction}

This paper is an expository paper which is the sequel of
\cite{ConcuToAlgTopo}. We will come back only very succinctly on the
explanations given in this latter. A technical appendix explains some
of the notions used in the core of the paper and fixes some notations.
A reader who would need more information about algebraic topology or
homological algebra could refer to \cite{May,Weibel,Rotman,hatcher}. A
reader who would need more information about the geometric point of
view of concurrency theory could refer to \cite{HDA,HDA2}.

The purpose is indeed to explain with much more details~\footnote{Even if
the limited  required number of pages for this paper too entails to make
some shortcuts.}
the speculations of the last paragraph of \cite{ConcuToAlgTopo}. More
precisely, we are going to describe a research program whose goal is
to construct an appropriate algebraic theory of the deformations of
higher dimensional automata (HDA) leaving invariant their
computer-scientific properties. Most of the  paper is as informal  as the
preceding one. The term \textit{dihomotopy} (contraction of
\textit{directed homotopy}) will be used as an analogue in our context
of the usual notion of \textit{homotopy}.

There are two known ways of modeling higher dimensional
automata for us to be able to study their deformations. 1) The
$\omega$-categorical approach, where strict globular
$\omega$-categories are supposed to encode the algebraic
structure of the possible compositions of execution paths and
homotopies between them, initiated by \cite{Pratt} and continued
in \cite{Gau} where connections with homological ideas of
\cite{HDA} were made. 2) The topological approach which consists,
loosely speaking, to locally endow a topological space with a
closed partial ordering which is supposed to represent the time :
this is the notion of local po-space developed for example in
\cite{HDA2}. The description of these models is sketched in
Section~\ref{formalization}.

Section~\ref{homological_construction} is an exposition of the
homological constructions which will play a role in the future
algebraic investigations.  Once again, the $\omega$-categorical case
and the topological case are described in parallel.

In Section~\ref{deform}, the notion of deformation of higher
dimensional automata is succinctly recalled. For further details,
see \cite{ConcuToAlgTopo}.

Afterwards Section~\ref{cat_ditype} exposes the main ideas about the
relation between homotopy types and dihomotopy types.  And some hints
to explore the algebraic structure of \textit{dihomotopy types} are
explained (this question is widely open).

Everything is presented in parallel because, as in usual algebraic
topology, the $\omega$-categorical approach and the topological
approach present a lot of similarities. In a first version, the paper
was organized with respect to the main result of \cite{KapVod2}, that
is the category equivalence between CW-complexes up to weak homotopy
equivalence and weak $\omega$-groupoids up to weak homotopy
equivalence.  By \cite{limit-groupoid}, it seems that this latter
result cannot be true, at least with the functors used in
Kapranov-Voevodsky's paper.  Therefore in this new version, the
presentation of some ideas is slightly changed. I thank Sjoerd Crans
for letting me know this fact.

\section{The formalization}\label{formalization}

\subsection{The  $\omega$-categorical approach}

Several authors have noticed that a higher dimensional automaton can
be encoded in a structure of precubical set
(Definition~\ref{def_cubique}). This idea is implemented in
\cite{cridlig96implementing} where a CaML program translating
programs written in Concurrent Pascal into (huge !) text files is
presented.

But such object does not contain any information about the way of
composing $n$-transitions, hence the idea of adding composition laws.
In an $\omega$-category (Definition~\ref{omega_categories}), the
$1$-morphisms represent the execution paths, the $2$-morphisms the
concurrent execution of the $1$-source and the $1$-target of the
$2$-morphisms we are considering, etc... The link between this way of
modeling higher dimensional automata and the formalization by
precubical sets is the realization functor $K\mapsto \Pi(K)$
described in Appendix~\ref{real}.

There exist two equivalent notions of (strict)
$\omega$-categories, the globular one and the cubical one
\cite{cub-glob} : the globular version will be used, although all
notions could be adapted to the cubical version. In fact, for
some technical reasons (Proposition~\ref{path-cat}), even a more
restrictive notion will be necessary\thinspace:

\bd\label{noncontractant}\cite{sglob} An $\omega$-category $\C$ is
\textit{non-contracting} if $s_1 x$ and $t_1 x$ are
$1$-dimen\-sional as soon as $x$ is not $0$-dimensional. Let $f$ be
an $\omega$-functor from $\C$ to $\D$.  The morphism $f$ is
\textit{non-contracting} if for any $1$-dimensional $x\in \C$, the
morphism $f(x)$ is a $1$-dimensional morphism of $\D$.
The category of non-contracting $\omega$-categories with the
non-contracting $\omega$-functors is denoted by $\omega Cat_1$.
\ed

The following proposition ensures that this technical  restriction is
not too small and that it does contain all precubical sets.

\bp\cite{sglob} For any precubical set $K$, $\Pi(K)$ is a non-contracting
$\omega$-category.  The functor $\Pi$ from the category of cubical
sets $Sets^{{\square^{pre}}^{op}}$ to that of $\omega$-catego\-ries $\omega Cat$
yields a functor from $Sets^{{\square^{pre}}^{op}}$ to the category of
non-contracting $\omega$-categories $\omega Cat_1$. \ep

\subsection{The topological approach}

Another way of modeling higher dimensional automata is to use the
notion of local po-space. A local po-space is a gluing of the
following local situation : 1) a topological space, 2) a partial
ordering, 3) as compatibility axiom between both structures, the
graph of the partial ordering is supposed to be closed
\cite{HDA2} (cf. Appendix~\ref{pospace}).

However the category of local po-spaces is too wide, and as in
usual algebraic topology, a more restrictive notion is necessary
to avoid too pathological situations (for instance think of the
Cantor set). A new notion which would play in this context the
role played by the CW-complexes in usual algebraic topology is
necessary. This is  precisely the subject of \cite{diCW}
(joined work with Eric Goubault).

Let $n\geq 1$. Let $D^n$ be the closed $n$-dimensional disk
defined by the set of points $(x_1,\dots,x_n)$ of $\R^n$ such
that $x_1^2+\dots +x_n^2\leq 1$ endowed with the topology induced
by that of $\R^n$. Let $S^{n-1}=\de D^n$ be the boundary of $D^n$
for $n\geq 1$, that is the set of $(x_1,\dots,x_n)\in D^n$ such
that $x_1^2+\dots +x_n^2=1$. Notice that $S^0$ is the discrete
two-point topological space $\{-1,+1\}$.  Let $I=[0,1]$. Let
$D^0$ be the one-point topological space. And let $e^n:=D^n-S^n$.
Loosely speaking, globular CW-complexes are gluing of po-spaces
$\vec{D}^{n+1}:=Glob(D^n)$ along $\vec{S}^{n}:=Glob(\de
D^n)=Glob(S^{n-1})$ where $Glob$ is the Globe functor (cf.
Appendix~\ref{pospace}).

Notice that there is a canonical inclusion of po-spaces
$\vec{S}^{n}\subset \vec{D}^{n+1}$ for $n\geq 1$. By convention, let
$\vec{S}^0:=\{0,1\}$ with the trivial ordering ($0$ and $1$ are not
comparable). There is a canonical inclusion  $\vec{S}^0\subset
\vec{D}^{1}$ which is a morphism of po-spaces.

\begin{propdef}\label{ndicell}\cite{diCW}
For any $n\geq 1$, $\vec{D}^{n}-\vec{S}^{n-1}$ with the induced
partial ordering is a po-space. It is called the $n$-dimensional
globular cell. More generally, every local po-space isomorphic to
$\vec{D}^{n}-\vec{S}^{n-1}$ for some $n$ will be called a
$n$-dimensional globular cell.
\end{propdef}

Now we are going to describe the process of attaching globular cells.

\begin{enumerate}
\item Start with a discrete set of points $X^0$.
\item Inductively, form the $n$-skeleton $X^n$ from $X^{n-1}$
by attaching globular $n$-cells $\vec{e}^n_\alpha$ via maps
$\phi_\alpha:\vec{S}^{n-1}\longrightarrow X^{n-1}$ with
$\phi_\alpha(\ei),\phi_\alpha(\ef)\in X^0$ such that\footnote{This
condition will appear to be necessary in the sequel.} : for every
non-decreasing map $\phi$ from $\vI$ to $\vec{S}^{n-1}$ such that
$\phi(0)=\ei$ and $\phi(1)=\ef$, there exists $0=t_0<\dots
<t_k=1$ such that $\phi_\alpha\circ \phi(t_i)\in X^0$ for any
$0\leq i\leq k$ which must satisfy
\begin{enumerate}
\item for any $0\leq i\leq k-1$, there exists a globular cell of
dimension $d_i$ with $d_i\leq n-1$ $\psi_i:\vec{D}^{d_i}\rightarrow
X^{n-1}$ such that for any $t\in[t_i,t_{i+1}]$, $\phi_\alpha\circ
\phi(t)\in \psi_i(\vec{D}^{d_i})$ ;
\item for $0\leq i\leq k-1$, the restriction of $\phi_\alpha\circ
\phi$ to $[t_i,t_{i+1}]$ is non-decreasing\thinspace;
\item the map $\phi_\alpha\circ \phi$ is non-constant ;
\end{enumerate}
Then $X^n$ is the quotient space of the disjoint union
$X^{n-1}\bigsqcup_{\alpha}\vec{D}^n_\alpha$ of $X^{n-1}$ with a
collection of $\vec{D}^n_\alpha$ under the identification $x\sim
\phi_\alpha(x)$ for $x\in \vec{S}^{n-1}_\alpha\subset
\de\vec{D}^n_\alpha$. Thus as set,
$X^n=X^{n-1}\bigsqcup_{\alpha}\vec{e}^n_\alpha$ where each
$\vec{e}^n_\alpha$ is a $n$-dimensional globular cell.
\item One can either stop this inductive process at a finite stage,
  setting $X=X^n$, or one can continue indefinitely, setting
  $X=\bigcup_n X^n$. In the latter case, $X$ is given the weak
  topology : A set $A\subset X$ is open (or closed) if and only if
  $A\cap X^n$ is open (or closed) in $X^n$ for some $n$ (this topology
  is nothing else but the direct limit of the topology of the $X^n$,
$n \in \N$). Such a $X$ is called a
  globular CW-complex and $X_0$ and the collection of
  $\vec{e}^n_\alpha$ and its attaching maps
  $\phi_\alpha:\vec{S}^{n-1}\longrightarrow X^{n-1}$ is called the
  cellular decomposition of $X$.
\end{enumerate}

\begin{figure}
\begin{center}
\includegraphics[width=7cm]{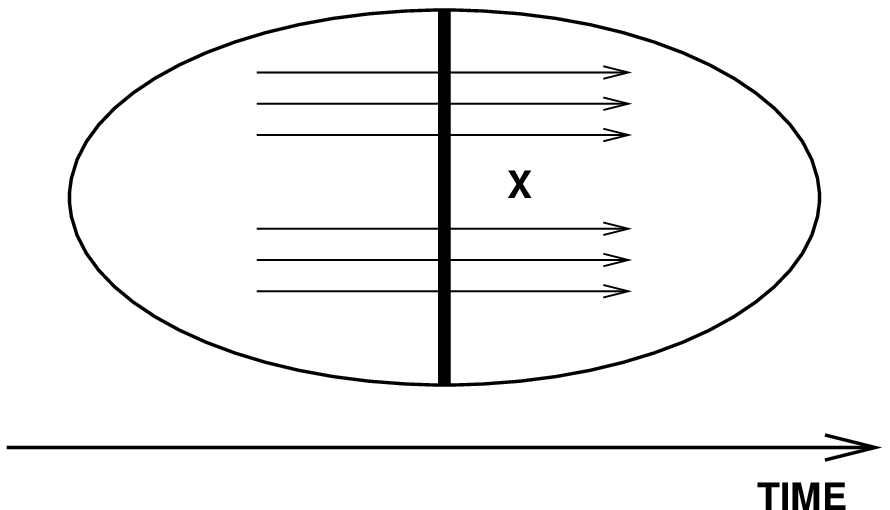}
\end{center}
\caption{Symbolic representation of $Glob(X)$ for some topological space $X$}
\label{2cell}
\end{figure}

As trivial examples of globular CW-complexes, there are
$\vec{D}^{n+1}$ and  $\vec{S}^{n}$ themselves where the
$0$-skeleton is, by convention, $\{\ei,\ef\}$.

We will consider without further mentioning that the segment $\vI$ is
a globular CW-complex, with $\{0,1\}$ as its $0$-skeleton.

\begin{propdef}\cite{diCW}
Let $X$ be a globular CW-complex with characteristic maps
$(\phi_\alpha)$. Let $\gamma$ be a continuous
map from $\vI$ to $X$. Then $\gamma([0,1])\cap X^0$ is finite.
Suppose that there exists $0\leq t_0<\dots<t_n\leq 1$ with
$n\geq 1$ such that $t_0=0$, $t_n=1$,  such that
for any $0\leq i\leq n$, $\gamma(t_i)\in X^0$, and at last
such that for any $0\leq i\leq n-1$,
there exists an $\alpha_i$ (necessarily unique) such that for $t\in [t_i,t_{i+1}]$,
$\gamma(t)\in \phi_{\alpha_i}(\vec{D}^{n_\alpha})$. Then such a
$\gamma$ is called an \textit{execution path} if
the restriction $\gamma\restriction_{[t_i,t_{i+1}]}$
is non-decreasing.
\end{propdef}

By constant execution paths, one means an execution paths $\gamma$ such that
$\gamma([0,1])=\{\gamma(0)\}$. The points (i.e. elements of the $0$-skeleton)
of a given globular CW-complexes $X$ are also called
\textit{states}. Some of them are fairly special:

\bd Let $X$ be a globular CW-complex. A point $\alpha$ of $X^0$ is
initial (resp. final) if for any execution path $\phi$ such that
$\phi(1)=\alpha$ (resp. $\phi(0)=\alpha$), then $\phi$ is the constant
path $\alpha$. \ed

Let us now describe the category of \textit{globular CW-complexes}.

\bd\cite{diCW} The category $\diCW$ of globular CW-complexes
is the category having as objects the globular CW-complexes and as
morphisms the continuous maps  $f:X\longrightarrow Y$
satisfying the
two following properties :
\begin{itemize}
\item $f(X^0)\subset Y^0$
\item for every non-constant execution path $\phi$ of $X$,
$f\circ \phi$ must not only be an execution path ($f$ must preserve
partial order), but also $f\circ \phi$ must be non-constant
as well  : we say that $f$ must be \textit{non-contracting}.
\end{itemize}
\ed

The condition of non-contractibility is very analogous to the notion
of non-contracting $\omega$-functors appearing in \cite{Gau}, and is
necessary for similar reasons. In particular, if the constant paths
are not removed from $\P^\pm X$ (see Section~\ref{hom-topo} for the
definition), then this latter spaces are homotopy equivalent to the
discrete set $X^0$ (the $0$-skeleton of $X$ !). And the removing of
the constant paths from $\P^\pm X$ entails to remove also the constant
paths from $\P X$ in order to keep the existence of both natural
transformations $\P\rightarrow \P^\pm$. Then the mappings $\P$ and
$\P^\pm$ can be made functorial only if we work with non-contracting
maps as above \cite{diCW}.

One can also notice that by construction, the attaching maps are
morphisms of globular CW-complexes. Of course one has

\bth\cite{diCW} Every globular CW-complex is a local po-space and this 
mapping induces a functor from the category of globular CW-complexes to 
the category of local po-spaces. \eth

\section{The homological constructions}\label{homological_construction}

The three principal constructions are all based upon the idea of
capturing the algebraic structure of the set of \textit{achronal
cuts} (cf. \cite{ConcuToAlgTopo} for some explanations of this
idea) included in the higher dimensional automaton $M$ we are
considering in three simplicial sets which seem to be the
basement of an algebraic theory which remains to build.  For
that, one has to construct in both cases (the categorical and the
topological approaches), three \textit{spaces} :
\begin{enumerate}
\item  the \textit{space of non-constant execution paths}
(this idea will become more precise below) : let us call it the \textit{path space} $\P M$
\item the  \textit{space of equivalence classes of non-constant execution paths beginning in the same way}
: let us call it the \textit{negative semi-path space} $\P^- M$
\item the  \textit{space of equivalence classes of non-constant execution paths ending in the same way}
: let us call it the \textit{positive semi-path space} $\P^+ M$.
\end{enumerate}
and one will consider the simplicial nerve of each one.

\begin{figure}
\[\xymatrix{& {\P M}\ar@{->}[ld]_{h^-}\ar@{->}[rd]^{h^+}&\\
{\P^- M}&&{\P^{+}M}}\]
\caption{The fundamental diagram}
\label{fundamental}
\end{figure}
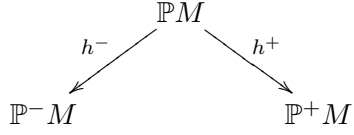

\cite{ConcuToAlgTopo} Figure~11  will become
Figure~\ref{fundamental} in both topological and
$\omega$-categorical situations. The construction of $h^-$ and
$h^+$ is  straightforward in both situations.

\subsection{The  $\omega$-categorical approach}

\bp\cite{sglob}\label{path-cat} Let $\C$ be an $\omega$-category.
Consider the set $\P\C=\bigcup_{n\geq 1} \C_n$.  Then the operators
$s_n$, $t_n$ and $*_n$ for $n\geq 1$ are internal to $\P\C$ if and
only if $\C$ is non-contracting. In that case, $\P\C$ can be endowed
with a structure of $\omega$-category whose $n$-source (resp.
$n$-target, $n$-dimensional composition law) is the $(n+1)$-source
(resp. $(n+1)$-target, $(n+1)$-dimensional composition law) of $\C$.
The $\omega$-category $\P\C$ is called the \textit{path
  $\omega$-category} of $\C$, and the mapping $\C$ induces a
well-defined functor from $\omega Cat_1$ to $\omega Cat$. \ep

\bd Let $\C$ be a non-contracting $\omega$-category. Denote by
$\mathcal{R}^-$ (resp. $\mathcal{R}^+$) the reflexive symmetric
and transitive closure of $\{(x,x *_0 y), x,y,x*_0y\in \P\C\}$
(resp. $\{(y,x *_0 y), x,y,x*_0y\in \P\C\}$) in $\P\C\p\P\C$. \ed

\bp\cite{fibrantcoin} Let $\alpha\in\{-,+\}$ and let $\C$ be a
non-contracting $\omega$-category. Then the universal problem
\begin{verse}
``There exists a pair $(\D,\mu)$ such that $\D$ is an
$\omega$-category and $\mu$ an $\omega$-functor from $\P\C$
to $\D$ such that for any $x,y\in \P\C$, $x\mathcal{R}^\alpha y$ implies
$\mu(x)=\mu(y)$.''
\end{verse} has a solution $(\P^\alpha\C,(-)^\alpha)$.
Moreover $\P^\alpha\C$ is generated by the elements of the form
$(x)^\alpha$ for $x$ running over $\P\C$. The mappings $\P^-$ and
$\P^+$ induce two well-defined functors from $\omega Cat_1$ to
$\omega Cat$. \ep

\bd The $\omega$-category $\P^-\C$ (resp. $\P^+\C$) is called the
negative (resp. positive) semi-path $\omega$-category of $\C$. \ed

In the sequel, $\P\C$ will be supposed to be a strict globular
$\omega$-groupoid in the sense of Brown-Higgins, which implies that 
$\P^-\C$ and $\P^+\C$ satisfy also the same property : this means concretely that 
if there exists an homotopy from a given execution path $\gamma$ to another one 
$\gamma'$, then there exists also an homotopy in the opposite direction
\cite{fibrantcoin}.


\bd\cite{sglob} The \textit{globular simplicial nerve}
$\mathcal{N}^{gl}$ is the functor from $\omega Cat_1$ to
$Sets^{\Delta^{op}}_+$ defined by $$\mathcal{N}^{gl}_n(\C):=\omega
Cat(\Delta^n,\P\C)$$ for $n\geq 0$ and with
$\mathcal{N}^{gl}_{-1}(\C):=\C_0\p\C_0$, and endowed with the
augmentation map $\de_{-1}$ from $\mathcal{N}^{gl}_0(\C)=\C_1$ to
$\mathcal{N}^{gl}_{-1}(\C)=\C_0\p \C_0$ defined by
$\de_{-1}x:=(s_0x,t_0 x)$. \ed

\begin{figure}
\begin{center}
\includegraphics[width=8cm]{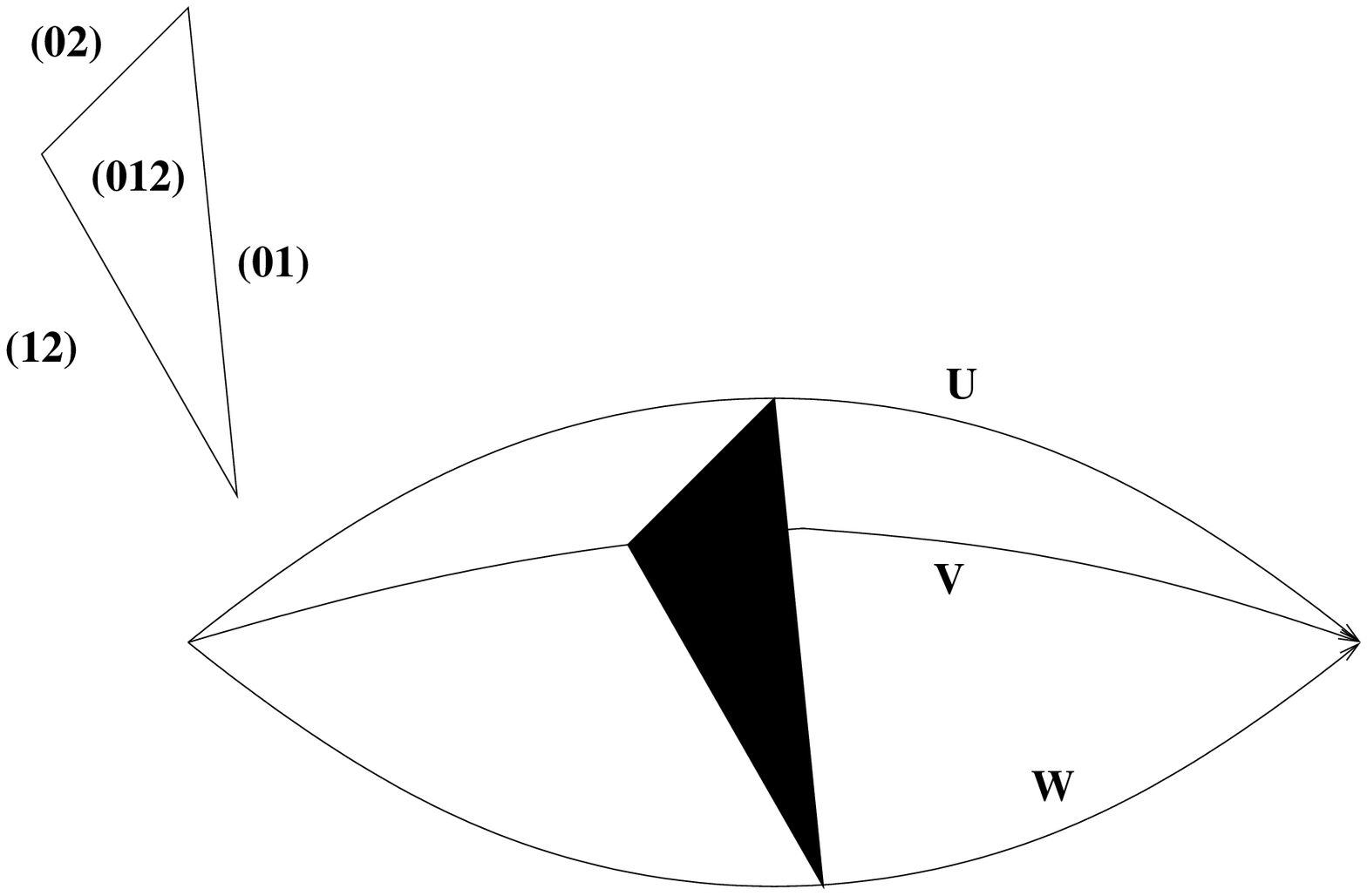}
\end{center}
\caption{Globular $2$-simplex} \label{2simplglob}
\end{figure}

Geometrically, a simplex of this simplicial nerve looks as in
Figure~\ref{2simplglob}.

\bd Let $\C$ be a non-contracting $\omega$-category. Then set
$$\mathcal{N}^{gl^-}_n(\C):=\omega Cat(\Delta^n,\P^-\C)$$
and $\mathcal{N}^{gl^-}_{-1}(\C):=\C_0$ with $\de_{-1}(x):=s_0 x$.
Then $\mathcal{N}^{gl^-}$ induces a functor from $\omega Cat_1$ to
$Sets^{\Delta^{op}}_+$ which is called the negative semi-globular nerve or (
branching semi-globular homology) of
$\C$ . \ed

\begin{figure}
\begin{center}
\includegraphics[width=10cm]{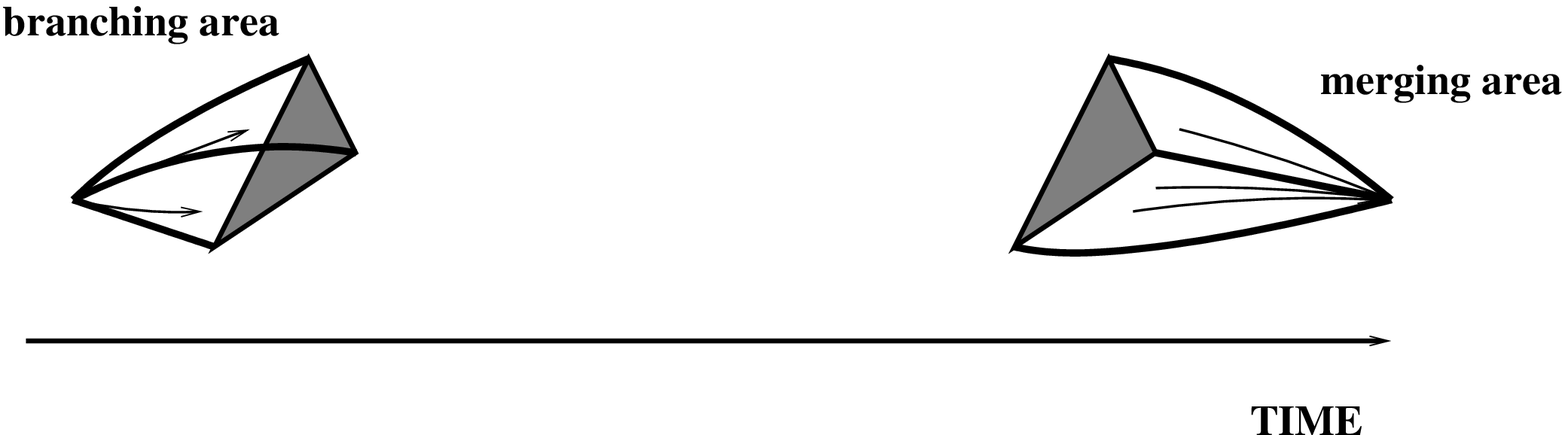}
\end{center}
\caption{Negative and positive semi-globular $2$-simplexes} \label{h}
\end{figure}

The \textit{positive semi-globular nerve} is defined in a similar way by replacing
$-$ by $+$ everywhere in the above definition and by setting
$\de_{-1}(x)=t_0 x$. Intuitively, the
simplexes in the semi-globular nerves look as in Figure~\ref{h} :
they correspond to the left or right half part of
Figure~\ref{2simplglob}.

\subsection{The topological approach}\label{hom-topo}

Let $X$ be a globular CW-complex. Let $\alpha,\beta\in X^0$. Denote by $(X,\alpha,\beta)^\bot$
the topological space of non-decreasing non-constant
continuous maps $\gamma$ from $[0,1]$ endowed with the usual order to $X$
such that $\gamma(0)=\alpha$ and $\gamma(1)=\beta$ and endowed with
the compact-open topology. Then

\bd Let $X$ be a globular CW-complex. Then the \textit{path space} of
$X$ is the disjoint union $$\P X=\bigsqcup_{(\alpha,\beta)\in X^0\p X^0}
(X,\alpha,\beta)^\bot$$ endowed with the disjoint union topology. \ed

Now denote by $(X,\alpha)^{\bot^-}$ (resp. $(X,\beta)^{\bot^+}$)
the topological space of non-decreasing non-constant continuous
maps from $[0,1]$ with the usual order to $X$ such that
$\gamma(0)=\alpha$ (resp. $\gamma(1)=\beta$), endowed with the
compact-open topology. Then

\bd Let $X$ be a globular CW-complex. Then the \textit{negative semi-path space} $\P^- X$
(resp. \textit{positive semi-path space} $\P^+ X$) of $X$  are defined by
\beas
&&\P^- X = \bigsqcup_{\alpha\in X^0}(X,\alpha)^{\bot^-}\\
&&\P^+ X = \bigsqcup_{\beta\in X^0} (X,\beta)^{\bot^+}
\eeas
endowed with the disjoint union topology. \ed

The reader can notice that in the topological context, we do not
need anymore to consider something like the equivalence relations
$\mathcal{R}^-$ and $\mathcal{R}^+$. The reason is that,
ideologically (``moralement'' in french !), a $1$-morphism is of
length $1$. On contrary, a non-constant execution path is
homotopic to any shorter execution path~\footnote{in a ``natural way'' by considering $H(\gamma(t),u)=\gamma(tu)$. It is the reason why $\P^-X$ and $\P^+X$ are 
homotopy equivalent to $X^0$ if one does not remove the constant paths 
from their definition.}


\bd\cite{sglob} The \textit{globular simplicial nerve}
$\mathcal{N}^{gl}$ is the functor from $\diCW$ to
$Sets^{\Delta^{op}}_+$ defined by
$$\mathcal{N}^{gl}_n(X):=S_n(\P X)$$ for $n\geq 0$ where
$S_*$ is the singular simplicial nerve (cf. Appendix~\ref{simplicial_explanation})
with
$\mathcal{N}^{gl}_{-1}(X):=X^0\p X^0$, and endowed with the
augmentation map $\de_{-1}$ from $\mathcal{N}^{gl}_0(X)=\P X$ to
$\mathcal{N}^{gl}_{-1}(X)=X^0\p X^0$ defined by
$\de_{-1}\gamma:=(\gamma(0),\gamma(1))$. \ed

\bd Let $X$ be a globular CW-complex. Then set
$$\mathcal{N}^{gl^-}_n(X):=S_n(\P^- X)$$ for $n\geq 0$
and $\mathcal{N}^{gl^-}_{-1}(X):=X^0$ with $\de_{-1}(\gamma):=\gamma(0)$.
Then $\mathcal{N}^{gl^-}$ induces a functor from $\diCW$ to
$Sets^{\Delta^{op}}_+$ which is called the branching semi-globular nerve of
$X$ . \ed

The \textit{merging semi-globular nerve} is defined in a similar way by replacing
$-$ by $+$ everywhere in the above definition and by setting
$\de_{-1}(\gamma):=\gamma(1)$.

\section{Deforming higher dimensional automata}\label{deform}

As already seen in \cite{ConcuToAlgTopo} in the
$\omega$-categorical context, there are two types of deformation
leaving invariant the computer scientific properties of higher
dimensional automata : the \textit{temporal deformations} (or
\textit{T-deformations}) and the \textit{spatial deformations}
(\textit{S-deformations}). The first type (temporal) is closely
related to the notion of homeomorphism because a non-trivial
execution path cannot be contracted in the same dihomotopy class~\footnote{In fact, 
the T-dihomotopy equivalences in \cite{diCW} are precisely the morphisms of globular CW-complexes 
inducing an homeomorphim between both underlying topological spaces.},
and the second one (spatial) to the classical notion of homotopy
equivalence.

\begin{figure}
\begin{center}
\subfigure[$\C$]{\label{dilate1}
\xymatrix{\alpha\fr{u}& \beta \ar@/^20pt/[r]^v \ar@/_20pt/[r]_w&\gamma}}
\hspace{2cm}
\subfigure[Subdivision  of $u$ in $\C$]{
\label{dilate2}
\xymatrix{\alpha_1\fr{u_1}&\alpha_2 \fr{u_2}& \beta \ar@/^20pt/[r]^v \ar@/_20pt/[r]_w&\gamma }}
\end{center}
\caption{Example of $T$-deformation}
\label{T1}
\end{figure}
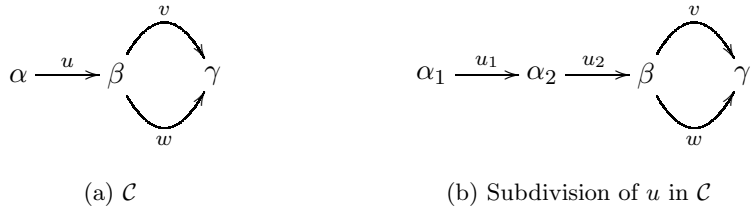

The $\omega$-categorical case will be only briefly recalled. A
temporal deformation corresponds informally to the reflexive
symmetric and transitive closure of subdividing in an
$\omega$-category a $1$-morphism in two parts as in Figure~\ref{T1}. A
spatial deformation consists of deforming in the considered
$\omega$-category $p$-morphisms with $p\geq 2$, which is equivalent to
deforming faces in one the three nerves in the usual sense of homotopy
equivalence.

The topological approach is completely similar. A temporal
deformation of a globular CW-complex $X$ consists of dividing in
two parts a globular $1$-dimensional cell of the cellular
decomposition of $X$, as in Figure~\ref{T1}. A spatial
deformation consists of crushing globular cells of higher
dimension.

Now what can we do with the previous homological constructions ?
First of all consider the corresponding simplicial homology
theories of all these augmented simplicial sets, with the
following convention on indices : for $n\geq -1$ and
$u\in\{gl,gl^-,gl^+\}$, set $H_{n+1}^u(M)=H_n\left(
  \mathcal{N}^u(M)\right)$ for $M$ either an $\omega$-category or a
globular CW-complex.  We obtain this way three homology theories
called as the corresponding nerve. One knows that the globular
homology sees the globes included in the HDA \cite{Gau,sglob} and
that the branching (resp. merging) semi-globular homology sees the branching
areas (resp. merging areas) in the HDA
\cite{Gau,Coin,fibrantcoin}. Since the three nerves are
Kan~\footnote{The $\omega$-categorical versions are Kan as soon
as $\P\C$ is an $\omega$-groupoid \cite{fibrantcoin} and the
singular simplicial nerve is known to be Kan \cite{May}.}, one can
also consider the homotopy groups of these nerves, with the same
convention for indices : for $n\geq 1$ and $u\in\{gl,gl^-,gl^+\}$,
set $\pi_{n+1}^u(X)=\pi_n\left(\mathcal{N}^u(X),\phi\right)$ for
$X$ either an $\omega$-category or a globular CW-complex. In this
latter case, the base-point $\phi$ is in fact a $0$-morphism of
$\P\C$, that is a $1$-morphism of $\C$ if $u=gl$, and an
equivalence class of $1$-morphism of $\C$ with respect to
$\mathcal{R}^-$ (resp. $\mathcal{R}^+$) if $u=gl^-$ (resp.
$u=gl^+$). Intuitively, elements of $\pi_{n+1}^{gl}$ are
$(n+1)$-dimensional cylinders with achronal basis.

The four first lines of Table~\ref{synthese} are explained in
\cite{ConcuToAlgTopo}.  The branching and merging (semi-cubical)
nerves $\mathcal{N}^\pm$ defined in \cite{Gau,Coin} 
are almost never Kan : in fact as soon as
there exists in the $\omega$-category $\C$ we are considering two
$1$-morphisms $x$ and $y$ such that $x *_0 y$ exists (see
Proposition~\ref{SurImageGlob}), both semi-cubical nerves are not Kan.

If the branching and merging (semi-cubical) nerves are replaced by the
branching and merging semi-globular nerve, then the ``almost'' (in
fact a ``no'') becomes a ``yes'' here because we are not disturbed
anymore by the non-simplicial part of the elements of the branching
and merging nerves (which is removed by construction).

The lines concerning the (globular, negative and positive semi-globular)
homotopy groups need to be explained. The S-invariance of a given
nerve  implies of course the S-invariance of the corresponding
homotopy groups.  As for the T-invariance, it is due to the fact
that in these homotopy groups, the ``base-point'' is an execution
path (or eventually an equivalence class of). So these homotopy
groups contain information only related to achronal cuts crossing
the ``base-point''. Dividing this base-point or any other
$1$-morphism or $1$-dimensional globular cell changes nothing.

The last lines are concerned with the bisimplicial set what we call
\textit{biglobular nerve} (for the contraction of bisimplicial
globular nerve) described in \cite{ConcuToAlgTopo,sglob} (cf.
Appendix~\ref{biglob}) and constructed by considering the structure of
augmented simplicial object of the category of small categories of the
globular nerve. The biglobular nerve inherits the S-invariance of the
globular nerve. And its T-invariance is due to the T-invariance of the
simplicial nerve functor of small categories.  The answer ``yes ?''
means that it is expected to find ``yes'' in some sense... It is worth
noticing that in a true higher dimensional automaton, $1$-morphisms
are never invertible because the time is not reversible. So one cannot
expect to find a Kan bisimplicial set in the usual sense of the
notion.

\begin{figure}
\begin{center}
\begin{tabular}{|l||c|c|c|c|}
\hline
Functors  & S-invariant & T-invariant & Kan & $-\infty \rightarrow +\infty$\\
\hline
$\mathcal{N}^{gl}$ &  yes & no & yes & yes \\
\hline
$\mathcal{N}^{\pm}$ & almost  & no & almost never & yes \\
\hline
$H^{gl}$ & yes & no &   no meaning & yes \\
\hline
$H^{\pm}$ & yes & yes &   no meaning & yes \\
\hline
$\mathcal{N}^{gl^\pm}$ & yes  & no & yes & yes \\
\hline
$H^{gl^\pm}$ & yes & yes & no meaning & yes \\
\hline
$\pi^{gl}$ & yes & yes & no meaning & no \\
\hline
$\pi^{gl^\pm}$ & yes & yes & no meaning & no \\
\hline
$\mathcal{N}^{bigl}$ & yes  & yes  & yes ? & yes \\
\hline
$H^{bigl}$ & yes  & yes  & no meaning & yes \\
\hline
$(\pi^{bigl})$ & (yes)  & (yes)  & (no meaning) & (yes)\\
\hline
\end{tabular}
\end{center}
\caption{Behavior w.r.t the two types of deformations}
\label{synthese}
\end{figure}

The last column is not directly concerned with the different types of
deformations of HDA, but rather by the question of knowing if the
functors contain information from $t=-\infty$ to $t=\infty$. The answer 
is yes everywhere except for the three homotopy groups functors : the latter 
contain indeed information only related to achronal cuts crossing the 
``base-point''. One can by the way notice that, in the $\omega$-categorical 
case :

\bp\label{composition_cylindre} 
Let $\phi$ and $\psi$ be two $1$-morphisms of a non-contracting 
$\omega$-category $\C$. Suppose that $\phi *_0 \psi$ exists. Then 
\begin{enumerate}
\item If $\phi *_0 \psi$ is $1$-dimensional, then the mapping 
$(x,y)\mapsto x*_0 y$ partially defined on $\C_n \p \C_n$ induces 
a morphism of groups 
$\pi_{n+1}^{gl}(\C,\phi)\p \pi_{n+1}^{gl}(\C,\psi)\rightarrow 
\pi_{n+1}^{gl}(\C,\phi *_0 \psi)$.
\item If $\phi *_0 \psi$ is $0$-dimensional, then the mapping 
$(x,y)\mapsto x*_0 y$ partially defined on $\C_n \p \C_n$ induces the 
constant map $\phi *_0 \psi$.
\end{enumerate}
\ep

\bpf It is due to the fact that 
\[\pi_{n+1}^{gl}(\C,\phi)=\{x\in \C_{n+1}, s_1x=s_2x=\dots=s_nx=t_1x=t_2x=\dots=t_nx=\phi\}\]
with $*_n$ for group law.
\epf

The above proposition is a hint to correct the drawbacks of the globular 
and semi-globular homotopy groups.

The last line $\pi^{bigl}$ is explained with Philosophy~\ref{diHu}.

\section{The category of dihomotopy types}\label{cat_ditype}

\subsection{Towards a construction}

Here both  approaches slightly diverge because of a lack of
knowledge about the $\omega$-categorical ways of  constructing homotopy
types. However one can certainly define in both contexts a notion of
\textit{weak dihomotopy equivalence} : see \cite{diCW} for the
topological context. Then let
\begin{itemize}
\item $\omega Grp$ be the category of strict globular $\omega$-groupoids
 with the $\omega$-functors as morphisms, and
$\Ho(\omega Grp)$ its localization by the weak homotopy equivalences
\item $\omega Cat_1^{Kan}$ be the category of non-contracting
  $\omega$-categories $\C$ such that $\P\C$ is an
  $\omega$-groupoid with the
  non-contracting $\omega$-functors as morphisms, and $\Ho(\omega Cat_1^{Kan})$ its
localization by the weak dihomotopy equivalences
\item $\CW$ the category of CW-complexes with the continuous maps as morphisms, and
$\Ho(\CW)$ its localization by the weak homotopy equivalences
\item $\diCW$ the category of globular CW-complexes with the morphisms
  of globular CW-complexes as morphisms, and $\Ho(\diCW)$ its
  localization by the weak dihomotopy equivalences.
\end{itemize}

\begin{philo} Both localizations $\Ho(\omega Cat_1^{Kan})$ and
  $\Ho(\diCW)$ contain the precubical sets modulo spatial and
  temporal deformations. However, due to the fact that strict globular
  $\omega$-groupoids do not represent all homotopy types
  \cite{Brown-Higgins0}, but only those having a trivial Whitehead
  product, $\Ho(\omega Cat_1^{Kan})$ could be  not big enough to
  construct an appropriate algebraic setting.
\end{philo}

After \cite{limit-groupoid}, it is clear that the $\omega$-categorical
realization functor described in Section~\ref{real} loses some
homotopical information and that keeping the complete information
requires to work with $\omega$-categories where the associativity of
$*_n$ is weakened for any $n\geq 1$. However, this lost homotopical
information is only related to the geometric situation in achronal
cuts. In particular, this realization functor does not contract
$1$-morphisms. Therefore the $\Ho(\omega Cat_1^{Kan})$ framework could
be sufficient to study questions concerning deadlocks or other similar
$1$-dimensional phenomena.

\bd The category $\Ho(\diCW)$ is called the category of dihomotopy
types. \ed

To describe the relation between the usual situation and the directed
situation, we need two last propositions and definitions :

\bp\cite{diCW} Let $X$ be a CW-complex. Let $Glob(X)^0=\{\ei,\ef\}$
where $\ei$ (resp. $\ef$) is the equivalence class of $(x,0)$ (resp.
$(x,1)$). Then the cellular decomposition of $X$ yields a cellular decomposition
of $Glob(X)$ and this way, $Glob(-)$ induces a functor from $\CW$ to $\diCW$. \ep

\bp Let $G$ be an object of $\omega Grp$. Then there exists a unique object $Glob(G)$ of
$\omega Cat_1^{Kan}$ such that $\P Glob(G)=G$, $Glob(G)_0=\{\alpha,\beta\}$ is a two-element
set, and such that
$s_0(Glob(G)\backslash \{\beta\})=\{\alpha\}$ and $t_0(Glob(G)\backslash \{\alpha\})=\{\beta\}$.
Moreover the mapping $Glob$ induces a functor from $\omega Grp$ to $\omega Cat_1^{Kan}$. \ep

Both $Glob$ functors (called \textit{Globe functors}) yield
two functors
$$\Ho(\CW)\rightarrow \Ho(\diCW)$$
and
$$\Ho(\omega Grp)\rightarrow
\Ho(\omega Cat_1^{Kan})$$

In the topological context, one has :

\bp \label{plonge}
\cite{diCW} Let $X$ and $Y$ be two CW-complexes. Then $X$ and
$Y$ are homotopy equivalent if and only if $Glob(X)$ and
$Glob(Y)$ are dihomotopy equivalent. Therefore the functor
$\Ho(\CW)\rightarrow \Ho(\diCW)$ is an embedding.
\ep

\begin{question} Is it possible to find an $\omega$-categorical construction
of $\Ho(\diCW)$ ? \end{question}

\subsection{Investigating the algebraic structure of the category of dihomotopy types}

One can check that in both topological and $\omega$-categorical situations, the following
fact holds

\bp\label{SurImageGlob} 
(partially in \cite{Coin}) Let $\alpha\in\{-,+\}$. The
morphism $h^\alpha$ induces an isomorphism of simplicial sets
(not of augmented simplicial sets for trivial reason !)
$\mathcal{N}^{gl}(Glob(M))\iso \mathcal{N}^{gl^\alpha}(Glob(M))$.
Moreover in the $\omega$-categorical case,
$\mathcal{N}^{gl}(Glob(M))\iso \mathcal{N}^{\alpha}(Glob(M))$
where $\mathcal{N}^{\alpha}$ are the branching or the merging
nerves (depending on the value of $\alpha$) of an
$\omega$-category as defined and studied in \cite{Gau,Coin}.
Moreover, this common simplicial set is homotopy equivalent to
the simplicial nerve of $M$.  \ep

This important proposition together with Proposition~\ref{plonge}
suggests us a way of investigating the
algebraic structure of the category of dihomotopy types.

\begin{philo}\label{philo1}
Let $\mathbf{Th}$ be a theorem (or a notion) in usual algebraic
topology,
  i.e. concerning the category of homotopy types.  Let
  $\mathbf{Th}^{di}$ be its lifting (i.e. its analogue) on the category of
  dihomotopy types. Then the statement $\mathbf{Th}^{di}$ must
  specialize into $\mathbf{Th}$ on the image of the Globe functor.
\end{philo}

Following Baues's philosophy \cite{homotopieabstraite}, a first goal
would be then to lift from the usual situation to the directed
situation the Whitehead theorem and the Hurewicz theorems.  Concerning
the last one, it would be first necessary to understand what is the
analogue~\footnote{The solution given in \cite{Gau} is naturally
  wrong : the morphisms $h^-$ and $h^+$ are not the analogues of the Hurewicz morphism.
    When \cite{Gau} was being written, It was not known that
  the correct definition of the globular homology would come from the
  simplicial homology of a simplicial nerve. Moreover the role of achronal cuts
  was also not yet understood. The globular homology was introduced as
  an answer of Goubault's suggestion of finding the analogue of the
  Hurewicz morphism in ``directed homotopy'' theory. Then starting from the
  principle that the branching and merging homology theories could be
  an analogue of the singular homology, I wondered whether it was
possible to construct  a morphism abutting to both corner homologies. The
  globular homology was then  designed to be  the source of this
  morphism.}  of the Hurewicz morphism for the category of dihomotopy
types.

\begin{philo}\label{diHu} The target of the Hurewicz morphism in the directed situation
  is likely to be the biglobular homology $H^{bigl}$. This new Hurewicz
  morphism must contain is some way all usual Hurewicz morphisms of
  all achronal cuts. At last, the source (let us denote it by
  $\pi^{bigl}$) of the Hurewicz morphism must be S-invariant,
  T-invariant and must contain information concerning the geometry of
  the HDA from $t=-\infty$ to $t=+\infty$.
\end{philo}

Suppose $n\geq 2$. After Proposition~\ref{composition_cylindre}, a
possible idea in the $\omega$-categorical case would be then to build
a chain complex of abelian groups by considering elements
\[(x_1,\dots,x_p)\in \pi_{n+1}^{gl}(\C,\phi_1)\p \dots \p  \pi_{n+1}^{gl}(\C,\phi_p)\]
for all $p$ and all $p$-uples $(\phi_1,\dots,\phi_p)$ such that
$\phi_1 *_0 \dots *_0 \phi_p$ exists and by considering the simplicial
differential map induced by $*_0$. Let us call the corresponding homology
theory the \textit{toroidal homology} $H^{tor}_{*}(\C)$. Of course
this construction makes sense only for $n\geq 2$ because the
$\pi_2^{gl}$ are not necessarily abelian. Then the classical Hurewicz
morphism induces a natural transformation from $H^{tor}_{*}$ to the
$E^2_{*,n+1}$-term of one of the canonical spectral sequences converging to $H^{bigl}$.

As explained in the introduction, the goal would be to reach an
homological understanding of the geometry of flows modulo
deformations. In particular, we would like to find exact sequences. It
is then reasonable to think that

\begin{philo}
An exact sequence $F_1(M) \rightarrow F_2(M) \rightarrow
  F_3(M)$ telling us something about flows $M$ of execution paths
  modulo spatial and temporal deformations must use functors $F_1$,
  $F_2$ and $F_3$ invariant by spatial and temporal deformations.
\end{philo}

The weakness of the internal structure of the globular nerve (it
is a disjoint union of simplicial sets), its non-invariance with
respect to temporal deformations, and its natural correction by
considering the \textit{biglobular nerve} suggests that the
\textit{biglobular homology} (the total homology of this
bisimplicial set) has more interesting homological properties than
the globular homology.

Concerning the biglobular nerve, it is worth noticing that this object
contains the whole information about the position of achronal
simplexes and about the temporal structure of the underlying higher
dimensional automaton. So in some sense, the biglobular nerve contains
everything related to the geometry of HDA. Since the biglobular nerve
is expected to be S-invariant and T-invariant, then it is natural to
ask the following question :

\begin{question}Is it possible to recover all other S-invariant and T-invariant
  functors from the biglobular nerve ? For example, is it possible to
  recover the semi-globular  homology theories ? \end{question}

Another natural question would be to relate a given dihomotopy type to
the underlying homotopy type (when the flow of execution paths is
removed). If the biglobular nerve really contains the complete
information, then it should be possible to recover from it the
underlying homotopy type.

\begin{figure}
\begin{center}
\includegraphics[width=6cm]{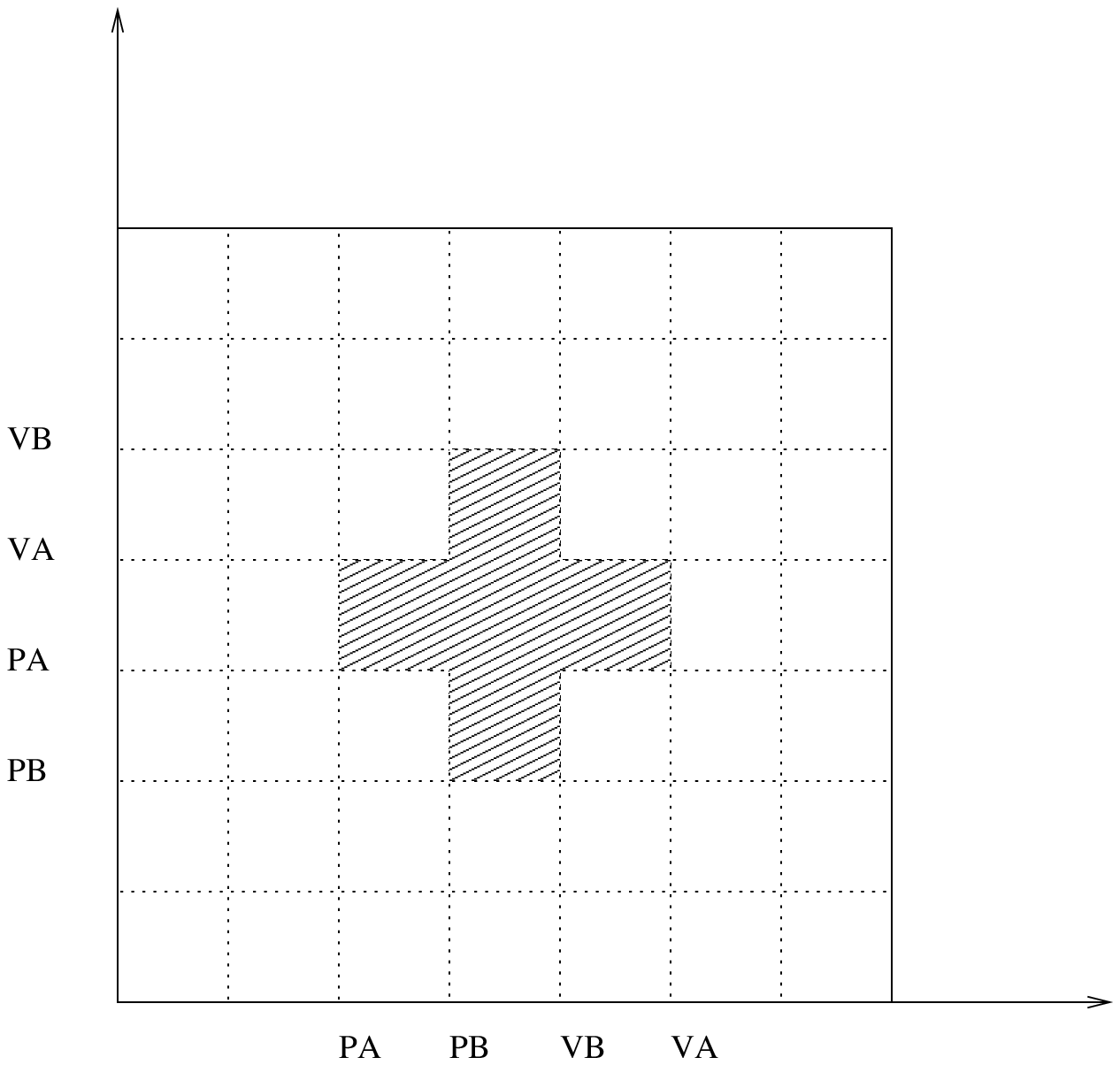}
\end{center}
\caption{PV diagram} \label{PV}
\end{figure}

As last remark, let us have a look at PV diagrams as in
Figure~\ref{PV}. They are always constructed by considering a
$n$-cube and by digging cubical holes inside.  Such examples produce
examples of $\omega$-categories or globular CW-complexes whose
all types of globular homologies do not have any torsion. To
classify this kind of examples, the study of \textit{rational
dihomotopy types} could be sufficient.

\section{Conclusion}

We have described in this paper a way of constructing the
category of dihomotopy types and we have given some hints to
investigate its internal algebraic structure. Intuitively, the
isomorphism classes of objects in this category represent exactly
the higher dimensional automata modulo the deformations which
leave invariant their computer-scientific properties. So a good
knowledge of the algebraic structure of this category will enable
us to classify higher dimensional automata up to dihomotopy and
therefore, hopefully, to write new algorithms manipulating
directly the equivalence classes of HDA.

\newpage

\appendix

\begin{center}
\textbf{\huge Technical Appendix}
\end{center}

\section{Local po-space : definition and examples}\label{pospace}

If $X$ is a topological space, a binary relation $R$ on $X$ is closed
if the graph of $R$ is a closed subset of the cartesian product $X\p
X$. If $R$ is a closed partial order $\leq$, then $(X,\leq)$ is called
a \textit{po-space} (see for instance \cite{nachbin}, \cite{johnstone}
and \cite{LFEGMRAlgebraic}).  Notice that a po-space is necessarily
Hausdorff. We say that $(U,\leq)$ is a sub-po-space of $(X,R)$ if and
only if it is a po-space such that $U$ is a sub topological space of
$X$ and such that $\leq$ is the restriction of $R$ to $U$.

A collection $\mathcal{U}(X)$ of po-spaces $(U,\leq_U)$ covering $X$ is called a
\textit{local partial order} if for every $x\in X$, there exists a
po-space $(W(x),\leq_{W(x)})$ such that:
\begin{itemize}
\item
 $W(x)$ is an open neighborhood
containing $x$,
\item the restrictions of $\leq_U$ and
$\leq_{W(x)}$ to $W(x)\cap U$ coincide for all $U\in
\mathcal{U}(X)$ such that $x\in U$. This can be stated as:
$y\leq_U z$ iff $y\leq_{W(x)} z$ for all $U\in \mathcal{U}(X)$
such that $x\in U$ and for all $y,z\in W(x)\cap U$.  Sometimes,
$W(x)$ will be denoted by $W_X(x)$ to avoid ambiguities. Such a
$W_X(x)$ is called a po-neighborhood.
\end{itemize}

Two local partial orders are equivalent if their union is a local
partial order. This defines an equivalence relation on the set of
local partial orders of $X$. A topological space together with an
equivalence class of local partial order is called a \textit{local
po-space}.

A morphism $f$ of local po-spaces (or \textit{dimap}) from
$(X,\mathcal{U})$ to $(Y,\Vt)$ is a continuous map from $X$ to
$Y$ such that for every $x\in X$,
\begin{itemize}
\item
there is a po-neighborhood $W(f(x))$ of $f(x)$ in $Y$,
\item
there exists a po-neighborhood $W(x)$ of $x$ in $X$ with
$W(x)\subset f^{-1}(W(f(x))$,
\item for $y,z\in W(x)$, $y\leq z$ implies
$f(y)\leq f(z)$.
\end{itemize}

In particular, a dimap $f$ from a po-space $X$ to a
po-space $Y$ is a continuous map from $X$ to $Y$ such that for any
$y,z\in X$, $y\leq z$ implies $f(y)\leq f(z)$. A morphism $f$
of local po-spaces from $[0,1]$ endowed with the usual ordering
(denoted by $\vec{I}$)
to a local po-space $X$ is called \textit{dipath} or sometime
\textit{execution path}.

The category of Hausdorff topological spaces with the continuous
maps as morphisms will be denoted by $\haus$. The category of
local po-spaces with the dimaps as morphisms will be denoted by
$\lpohaus$. The category of general topological spaces without
further assumption will be denoted by $\top$ and the category of
general topological spaces endowed with a partial ordering not
necessary closed will be denoted by $\potop$.

We end this section by an example of po-spaces which matters for this
paper. Let us construct the \textit{Globe} $Glob(X)$ associated to a
topological space $X$. It is defined as follows. As topological space,
$Glob(X)$ is the quotient of the product space $X\p I$ by the
relations $(x,0)=(x',0)$ and $(x,1)=(x',1)$ for any $x,x'\in X$. It is
equipped with the closed partial order $(x,t)\leq (x',t')$ if and only
if $x=x'$ and $t\leq t'$. The equivalence class of $(x,0)$ (resp.
$(x,1)$) in $Glob(X)$ is denoted by $\ei$ (resp. $\ef$).

\section{Simplicial set}\label{simplicial_explanation}

For further details, cf. \cite{May,Weibel}.

\bd A simplicial set $A_*$  is a family $(A_n)_{n\geq 0}$ together with
face maps $\de_i:A_n\rightarrow A_{n-1}$ and
$\epsilon_i:A_n\rightarrow A_{n+1}$ for $i=0,\dots,n$ which satisfy
the following identities :
\begin{alignat*}{2}
& \de_i\de_j = \de_{j-1}\de_i && \hbox{ if }i<j\\
&\epsilon_i \epsilon_j= \epsilon_{j+1}\epsilon_i && \hbox{ if }i\leq j\\
& \de_i \epsilon_j = \left\{\begin{array}{c} \epsilon_{j-1}\de_i
\\ \hbox{Identity} \\ \epsilon_j \de_{i-1} \end{array}\right. &&
\begin{array}{l} \hbox{ if }i<j \\ \hbox{ if }i=j\hbox{ or }i=j+1\\
\hbox{ if } i>j+1\end{array}
\end{alignat*}
A morphism of simplicial sets from $A_*$ to $B_*$ consists of a set
map from $A_n$ to $B_n$ for each $n\geq 0$ commuting with all
operators defined on both sides. The category of simplicial sets is denoted
by $Sets^{\Delta^{op}}$.
\ed

Consider the \textbf{topological $n$-simplex} $\Delta^n$ defined by
$$\Delta^n=\{(t_0,\dots,t_n), t_0\geq 0,\dots,t_n\geq0\hbox{ and } t_0+\dots t_n=1\}$$
Here is now the most classical example of simplicial sets :

\bd Let $Y$ be a topological space. The singular simplicial nerve
of $Y$ is the simplicial set $S_*(XY)$ defined as follows :
$S_n(Y):=\top(\Delta^n,Y)$ with
$\de_i(f)(t_0,\dots,t_{n-1})=f(t_0,\dots,t_{i-1},\dots,t_{n-1})$
and
$\epsilon_i(f)(t_0,\dots,t_{n+1})=f(t_0,\dots,t_{i-1},t_i+t_{i+1},t_{i+2},\dots,t_{n+1})$.
\ed

\bd\cite{triple} An augmented simplicial set is a simplicial set
\[((X_n)_{n\geq 0}, (\de_i:X_{n+1}\longrightarrow X_n)_{0\leq i\leq n+1},
(\epsilon_i:X_{n}\longrightarrow X_{n+1})_{0\leq i\leq n})\]
together with an additional set $X_{-1}$ and an additional map
$\de_{-1}$ from $X_0$ to $X_{-1}$ such that
$\de_{-1}\de_0=\de_{-1}\de_1$. A morphism of augmented simplicial
set is a map of $\N$-graded sets which commutes with all face and
degeneracy maps. We denote by $Sets^{\Delta^{op}}_+$ the category
of augmented simplicial sets.  \ed

If $X_*$ is an augmented simplicial set, one obtains a chain
complex of abelian groups ($\Z S$ being the free abelian group
generated by the set $S$)
\[
\xymatrix{  \dots \fr{}& {\Z X_2}\frr{\de_0-\de_1+\de_2} &&{\Z X_1}\frr{\de_0-\de_1} && {\Z X_0}\fr{\de_{-1}}& {\Z X_{-1}}\fr{} & 0 }\]
We will denote $H_{n+1}(X)$ for $n\geq -1$ the $n$-th simplicial homology group of
$X_*$. This means for example that $H_1(X)$ will be the quotient
 of $\de_{-1}:\mathcal{N}^{gl}_0(\C)\rightarrow \mathcal{N}^{gl}_{-1}(\C)$
by the image of $\de_0-\de_1:\mathcal{N}^{gl}_1(\C)\rightarrow
\mathcal{N}^{gl}_{0}(\C)$.

\section{Precubical set, globular $\omega$-category and globular set}\label{main_def}

\bd\label{def_cubique}\cite{Brown_cube} \cite{cube} A
\textit{precubical set} consists of a family of sets
$(K_n)_{n\geqslant 0}$ and of a family of face maps
$\xymatrix@1{{K_n}\fr{\de_i^\alpha} &{K_{n-1}}}$ for
$\alpha\in\{-,+\}$ which satisfies the following axiom (called
sometime the cube axiom) :
\begin{center}
$\de_i^\alpha \de_j^\beta = \de_{j-1}^\beta \de_i^\alpha$ for all
$i<j\leqslant n$ and $\alpha,\beta\in\{-,+\}$.
\end{center}
\ed

If $K$ is a precubical set, the elements of $K_n$ are called the
$n$-cubes. An element of $K_n$ is of dimension $n$. The elements
of $K_0$ (resp. $K_1$) can be called the \textit{vertices} (resp.
the \textit{arrows}) of $K$.

\bd\label{omega_categories}
\cite{Brown-Higgins0,oriental,Tensor_product} An
\textit{$\omega$-category} is a set $A$ endowed with two families
of maps $(s_n=d_n^-)_{n\geqslant 0}$ and $(t_n=d_n^+)_{n\geqslant
0}$ from $A$ to $A$ and with a family of partially defined 2-ary
operations $(*_n)_{n\geqslant 0}$ where for any $n\geqslant 0$,
$*_n$ is a map from $\{(a,b)\in A\p A, t_n(a)=s_n(b)\}$ to $A$
($(a,b)$ being carried over $a *_n b$) which satisfies the
following axioms for all $\alpha$ and $\beta$ in $\{-,+\}$ :

\begin{enumerate}
\item $d_m^\beta d_n^\alpha x=
\left\{\begin{CD}d_m^\beta x \hbox{  if $m<n$}\\  d_n^\alpha x
\hbox{  if $m\geqslant n$}
\end{CD}\right.$
\item $s_n x *_n x= x *_n t_n x = x$
\item if  $x *_n y$ is well-defined, then  $s_n(x *_n y)=s_n x$, $t_n(x *_n y)=t_n y$
and for  $m\neq n$, $d_m^\alpha(x *_n y)=d_m^\alpha x *_n
d_m^\alpha y$
\item as soon as the two members of the following equality exist, then
$(x *_n y) *_n z= x *_n (y *_n z)$
\item if $m\neq n$ and if the two members of the equality make sense, then
$(x *_n y)*_m (z*_n w)=(x *_m z) *_n (y *_m w)$
\item for any  $x$ in $A$, there exists a natural number $n$ such that $s_n x=t_n x=x$
(the smallest of these numbers is called the dimension of $x$ and
is denoted by $dim(x)$).
\end{enumerate}
\ed

A $n$-dimensional element of $\C$  is called a $n$-morphism. A
$0$-morphism is also called a state of $\C$, and a $1$-morphism an
arrow. If $x$ is a morphism of an $\omega$-category $\C$, we call
$s_n(x)$ the $n$-source of $x$ and $t_n(x)$ the $n$-target of $x$.
The category of all $\omega$-categories (with the obvious
morphisms) is denoted by $\omega Cat$. The corresponding morphisms
are called $\omega$-functors. The set of $n$-dimensional morphisms
of $\C$ is denoted by $\C_n$.

\begin{figure}
\begin{center}
\subfigure[Composition of two
$2$-morphims]{\label{composition-of-2trans}
\xymatrix{&&\\
{\alpha} \rruppertwocell<10>{A} \rrlowertwocell<-10>{B} \ar[rr]&&
{\beta} }} \hspace{2cm} \subfigure[The  $\omega$-category
$\Delta^2$]{\label{2simpl} \xymatrix{&&
{(0)}\\&&\\{(2)}\ar@{->}[rruu]^{(02)}\ar@{->}[rr]_{(12)}&&{(1)}\ar@{->}[uu]_{(01)}
\ar@2{->}[lu]^{(012)}}} \subfigure[The $\omega$-category
$I^3$]{\label{I3} \xymatrix{ &\ar@{->}[rr]^{+0-}&
&\ar@{->}[rd]^{++0}&&&
&\ar@{->}[dr]|{+-0}\ar@{->}[rr]^{+0-}&&\ar@{->}[dr]^{++0}&\\
\ar@{->}[ru]^{0--}\ar@{->}[rr]|{-0-}\ar@{->}[rd]_{--0}&&\ar@{->}[ru]|{0+-}\ar@{->}[dr]|{-+0}\ff{lu}{00-}&&
\ar@3{->}[rr]^{000}&&
\ar@{->}[ur]^{0--}\ar@{->}[dr]_{--0}&&\ff{ur}{+00}\ar@{->}[rr]|{+0+}&&\\
&\ar@{->}[rr]_{-0+}\ff{ru}{-00}&&\ar@{->}[ru]_{0++}\ff{uu}{0+0}&&&
&\ar@{->}[ru]|{0-+}\ar@{->}[rr]_{-0+}\ff{uu}{0-0}&&\ff{ul}{00+}\ar@{->}[ru]_{0++}\\
}}
\end{center}
\caption{Some $\omega$-categories (a $k$-fold arrow symbolizes a
k-morphism)}
\end{figure}

As fundamental examples of $\omega$-categories, there is the
$\omega$-category $\Delta^n$ freely generated by the faces of the
$n$-simplex \cite{oriental}.  To characterize this
$\omega$-category, the first step consists of labeling all faces
of the $n$-simplex. Its faces are indeed in bijection with
strictly increasing sequences of elements of $\{0,1,\dots,n\}$. A
sequence of length $p+1$ will be of dimension $p$. If $x$ is a
face, let $R(x)$ be the set of faces of $x$ seen  as a
sub-simplex. If $X$ is a set of faces, then let
$R(X)=\bigcup_{x\in X}R(x)$. Notice that $R(X\cup Y)=R(X)\cup
R(Y)$ and that $R(\{x\})=R(x)$. Then $\Delta^n$ is  the free
$\omega$-categories generated by the $R(x)$ with the rules

\begin{enumerate}
\item  For $x$ $p$-dimensional with $p\geq 1$,
$s_{p-1}(R(x))=R(s_x)$ and $t_{p-1}(R(x))=R(t_x)$ where $s_x$ and
$t_x$ are the sets of faces defined below.
\item If $X$ and $Y$ are two elements of $\Delta^n$ such that
$t_p(X)=s_p(Y)$ for some $p$, then $X\cup Y$ belongs to $\Delta^n$
and $X\cup Y=X *_p Y$.
\end{enumerate}

Let us give the definition of $s_x$ and $t_x$ on some example :
$$s_{(04589)}=\{(4589),(0489),(0458)\}$$
The elements in odd position
are removed ;
$$t_{(04589)}=\{(0589),(0459)\}$$
The elements in even
position are removed.

Let $\underline{\Delta}$ be the unique small category such that a
pre-sheaf over $\underline{\Delta}$ is exactly a simplicial set
\cite{May,Weibel}. The category $\underline{\Delta}$ has for
objects the finite ordered sets $[n]=\{0<1<\dots<n\}$ for integers
$n\geq 0$ and has for morphisms the non-decreasing monotone
functions. One is  used to distinguishing in this category the
morphisms $\epsilon_i:[n-1]\rightarrow [n]$ and
$\eta_i:[n+1]\rightarrow [n]$ defined as follows for each $n$ and
$i=0,\dots,n$ :
\[
\epsilon_i(j)=\left\{\begin{array}{c}j \hfill\hbox{ if
}j<i\\j+1\hfill\hbox{ if }j\geq i\end{array}\right\} ,\ \ \
\eta_i(j)=\left\{\begin{array}{c}j \hfill\hbox{ if }j\leq
i\\j-1\hfill\hbox{ if }j> i\end{array}\right\}
\]

The mapping  $n\mapsto \Delta^n$ yields a functor from
$\underline{\Delta}$ to $\omega Cat$ by setting $\epsilon_i
\mapsto \Delta^{\epsilon_i}$ and $\eta_i \mapsto \Delta^{\eta_i}$
where
\begin{itemize}
\item for any face $(\sigma_0<\dots <\sigma_s)$ of $\Delta^{n-1}$,
$\Delta^{\epsilon_i}(\sigma_0<\dots <\sigma_s)$ is the only face
of $\Delta^{n}$ having $\epsilon_i\{\sigma_0,\dots,\sigma_s\}$ as
set of vertices ;
\item for any face $(\sigma_0<\dots <\sigma_r)$ of $\Delta^{n+1}$,
$\Delta^{\eta_i}(\sigma_0<\dots <\sigma_r)$ is the only face of
$\Delta^{n}$ having $\eta_i\{\sigma_0,\dots,\sigma_r\}$ as set of
vertices.
\end{itemize}

Therefore

\bd Let $\C$ be an $\omega$-category. Then the graded set $\omega
Cat(\Delta^*,\C)$ is naturally endowed with a structure of
simplicial sets. It is called the simplicial nerve of $\C$. \ed

\section{$\omega$-categorical realization of a precubical
set}\label{real}

Intuitively the $\omega$-categorical realization $\Pi(K)$ of a
precubical set $K$ (also called the free $\omega$-category
generated by $K$) as defined below contains as $n$-morphisms all
composites (or all concatenations) of cubes of $K$ which are
$n$-dimensional (this means that somewhere in the composite, a
$n$-dimensional cube appears). In particular the $1$-morphisms of
$\Pi(K)$ will be exactly all arrows of $K$ and all possible
compositions of these arrows.

The free $\omega$-category $\Pi(K)$ is constructed as follows. The
main ingredient is the free $\omega$-category $I^n$ generated by
the faces of the $n$-cube. Its characterization is very similar to
that of the $\omega$-category $\Delta^n$ generated by the faces of
the $n$-simplex. The faces of the $n$-cube are labeled by the
word of length $n$ in the alphabet $\{-,0,+\}$, the number of zero
corresponding to the dimension of the face. Everything is
similar, except the definition of $s_x$ and $t_x$. The set $s_x$
is the set of sub-faces of the faces obtained by replacing the
$i$-th zero of $x$ by $(-)^i$, and the set $t_x$ is the set of
sub-faces of the faces obtained by replacing the $i$-th zero of
$x$ by $(-)^{i+1}$. For example,
$s_{0+00}=\{\hbox{-+00},\hbox{0++0},\hbox{0+0-}\}$ and
$t_{0+00}=\{\hbox{++00},\hbox{0+-0},\hbox{0+0+}\}$.
Figure~\ref{I3} represents the free $\omega$-category generated
by the $3$-cube (cf. \cite{Gau} for some examples of
calculations). The first construction of $I^n$ is due to
Aitchison in \cite{explicite}.

Then to each $x\in K_n$, we associate a copy of $I^n$ denoted by
$\{x\}\p I^n$ whose  corresponding faces will be denoted by
$(x,k_1\dots k_n)$. We then take the quotient of the direct sum of
these $\{x\}\p I^{dim(x)}$ in $\omega Cat$ (which corresponds for
the underlying sets to the disjoint union) by the relations
$$(\de_i^\alpha(x),k_1\dots k_{n-1})\sim (x,k_1\dots k_{i-1}
[\alpha]_i k_i \dots k_{n-1})$$ for any $n\geq 1$ and any $x\in
K_n$ where the notation $[\alpha]_i$ means that $\alpha$ is put in
$i$-th position. This expression means that in the copy of
$I^{n-1}$ corresponding to $\de_i^\alpha(x)$, the face $k_1\dots
k_{n-1}$ must be identified to the face $k_1\dots
k_{i-1}[\alpha]_i k_i \dots k_{n-1}$ in the copy of $I^n$
corresponding to $x$. And one has

\bp One obtains a well-defined $\omega$-category $\Pi(K)$ and
$\Pi$ induces a well-defined functor from the category of
precubical sets to that of $\omega$-categories. \ep

The proof uses the coend construction  (cf. \cite{cat}).

\section{Localization of a category with respect to a collection of morphisms}\label{def_localization}

\bd Let $\C$ be a category (not necessarily small). Let $S$ be a collection of morphisms
of $\C$. Consider the following universal problem :
\begin{verse}
``There exists a pair $(\D,\mu)$ such that $\mu$ is a functor from $\C$ to $\D$ and such that
for any $s\in S$, $\mu(s)$ is an invertible morphism of $\D$.''
\end{verse}
Then the solution $(\C[S^{-1}],Q)$, if there exists, is called the
localization of $\C$ with respect to $S$. \ed

\section{The biglobular nerve}\label{biglob}

\bth\cite{sglob}\label{grading} Let $\C$ be a non-contracting  $\omega$-category.
\begin{enumerate}
\item\label{extremite} Let $x$ be an $\omega$-functor from $\Delta^n$ to $\P\C$ for some
$n\geq 0$. Then the set maps $(\sigma_0\dots \sigma_r)\mapsto s_0 x((\sigma_0\dots \sigma_r))$
and $(\sigma_0\dots \sigma_r)\mapsto t_0 x((\sigma_0\dots \sigma_r))$ from
the underlying set of faces of $\Delta^n$ to $\C_0$ are constant. The unique value of
$s_0\circ x$ is denoted by $S(x)$ and the unique value of $t_0 \circ x$ is denoted by
$T(x)$.
\item For any pair $(\alpha,\beta)$ of $0$-morphisms of $\C$, for any
$n\geq 1$, and for any $0\leq i\leq n$, then
$\de_i \left(\mathcal{N}_n^{gl}(\C[\alpha,\beta])\right)\subset \mathcal{N}_{n-1}^{gl}(\C[\alpha,\beta])$.
\item For any pair $(\alpha,\beta)$ of $0$-morphisms of $\C$, for any
$n\geq 0$, and for any $0\leq i\leq n$, then
$\epsilon_i\left(\mathcal{N}_n^{gl}(\C[\alpha,\beta])\right) \subset \mathcal{N}_{n+1}^{gl}(\C[\alpha,\beta])$.
\item By setting, $G^{\alpha,\beta}\mathcal{N}_n^{gl}(\C):=\mathcal{N}_n^{gl}(\C[\alpha,\beta])$
for $n\geq 0$ and $G^{\alpha,\beta}\mathcal{N}_{-1}^{gl}(\C):=\{(\alpha,\beta),(\beta,\alpha)\}$, one
obtains a $(\C_0\p \C_0)$-graduation on the globular nerve ; in particular, one has the direct sum
of augmented simplicial sets
\[\mathcal{N}_*^{gl}(\C)=\bigsqcup_{(\alpha,\beta)\in\C_0\p\C_0} G^{\alpha,\beta}\mathcal{N}_*^{gl}(\C)\]
and $G^{\alpha,\beta}\mathcal{N}_*^{gl}(\C)=\mathcal{N}_*^{gl}(\C[\alpha,\beta])$.
\end{enumerate}
\eth

Let $\C$ be a non-contracting $\omega$-category. Using
Theorem~\ref{grading}, recall that for some $\omega$-functor $x$ from $\Delta^n$ to
$\P\C$, one calls $S(x)$ the unique element of the image of
$s_0\circ x$ and $T(x)$ the unique element of the image of
$t_0\circ x$. If $(\alpha,\beta)$ is a pair of $\mathcal{N}_{-1}^{gl}(\C)$,
set $S(\alpha,\beta)=\alpha$ and $T(\alpha,\beta)=\beta$.

\bp\cite{sglob} Let $\C$ be a non-contracting $\omega$-category. Let $x$
and $y$ be two $\omega$-functors from $\Delta^n$ to $\P\C$ with $n\geq
0$. Suppose that $T(x)=S(y)$.
Let $x*y$ be the map from the faces of
$\Delta^n$ to $\C$ defined by \[(x*y)((\sigma_0\dots
\sigma_r)):=x((\sigma_0\dots \sigma_r))*_0 y((\sigma_0\dots
\sigma_r)).\]
Then the following conditions are equivalent :
\begin{enumerate}
\item\label{c1} The image of $x*y$ is a subset of $\P\C$.
\item\label{c2} The set map $x*y$ yields an $\omega$-functor from $\Delta^n$ to $\P\C$
and $\de_i(x*y)=\de_i(x)*\de_i(y)$ for any $0\leq i\leq n$.
\end{enumerate}
On contrary, if for some $(\sigma_0\dots \sigma_r)\in \Delta^n$,
$(x*y)((\sigma_0\dots \sigma_r))$ is $0$-dimensional, then $x*y$ is
the constant map $S(x)=T(y)$.
\ep

In the sequel, we set $(\alpha,\beta)*(\beta,\gamma)=(\alpha,\gamma)$,
$S(\alpha,\beta)=\alpha$ and $T(\alpha,\beta)=\beta$.
If $x$ is an $\omega$-functor
from $\Delta^n$ to $\P\C$, and if $y$ is the constant map $T(x)$ (resp. $S(x)$)
from $\Delta^n$ to $\C_0$, then set $x*y:=x$ (resp. $y*x:=x$).

\bth Suppose that $\C$ is an object of $\omega Cat_1$. Then for
$n\geq 0$, the operations $S$, $T$ and $*$ allow to define a
small category $\underline{\mathcal{N}^{gl}_n(\C)}$ whose
morphisms are the elements of $\mathcal{N}^{gl}_n(\C)\cup
\{\hbox{constant maps }\Delta^n\rightarrow \C_0\}$ and whose
objects are the $0$-morphisms of $\C$. If
$\underline{\mathcal{N}^{gl}_{-1}(\C)}$ is the small category
whose morphisms are the elements of $\C_0\p\C_0$ and whose
objects are the elements of $\C_0$ with the operations $S$, $T$
and $*$ above defined, then one obtains (by defining the face
maps $\de_i$ and degeneracy maps $\epsilon_i$ in an obvious way
on $\{\hbox{constant maps }\Delta^n\rightarrow \C_0\}$) this way
an augmented simplicial object $\underline{\mathcal{N}^{gl}_*}$
in the category of small categories. \eth

By composing by the classifying space functor of small categories
(cf. for example \cite{classifiant} for further details), one obtains
a bisimplicial set which is called the \textit{biglobular nerve}.

\end{document}